\documentclass[11pt]{article}

\usepackage{amsmath,amssymb,amsthm,overpic,subfig,booktabs,amsrefs,color,booktabs}

\def\Bmeas{\mathcal{B}}
\def\NNN{\mathcal{N}}
\def\bf0{\mathbf{0}}
\def\bfc{\mathbf{c}}
\def\bfX{\mathbf{X}}
\def\bfx{\mathbf{x}}
\def\bfY{\mathbf{Y}}
\def\bfy{\mathbf{y}}
\def\bfZ{\mathbf{Z}}
\def\bfz{\mathbf{z}}
\def\bfv{\mathbf{v}}
\def\bfQ{\mathbf{Q}}
\def\bfD{\mathbf{D}}
\def\bfmu{{\boldsymbol\mu}}
\def\bfSigma{{\boldsymbol\Sigma}}

\def\diffns{\mathrm{d}}
\def\diag{\mathrm{diag}}
\def\dd{\mathrm{d}}

\def\TT{\mathrm{T}}

\def\Pmeas{\mathbb{P}}
\def\RR{\mathbb{R}}

\begin{document}
\title{A Chi Distribution Model of Hail Storm Damage}
\date{}

\author{Paul C. Kettler\footnote{Centre of Mathematics for Applications, Department of Mathematics, University of Oslo, P.O. Box 1053, Blindern, N--0316 Oslo, Norway, \texttt{paulck\@@math.uio.no}} \ and Georg Muntingh\footnote{SINTEF ICT, P.O. Box 124 Blindern, 0314 Oslo, Norway, and Centre of Mathematics for Applications, Department of Mathematics, University of Oslo, P.O. Box 1053, Blindern, N--0316 Oslo, Norway, \texttt{georgmu\@@math.uio.no}}}

\maketitle

\begin{abstract}
\noindent This paper addresses the pattern of damage, and investigates its properties, of a theoretical hail storm which gathers in intensity before subsiding, and which travels linearly across the landscape at constant velocity.  We start by assuming a simpler model, that of a storm which does not move, restricted to having an uncorrelated binormal distribution of damage. This model, expressed in the natural polar coordinates, leads to a 1-dimensional pattern of damage as a function of the marginal radial distance conforming to the $\chi$-distribution with two degrees of freedom. The model is then extended to the traveling form, allowing further for a correlation of the variables, extending, as well, to the multidimensional case. In its full florescence the model produces hyperellipsoidal hypersurfaces of equal intensity for the correlated multinormal assumption. We provide closed-form solutions for the totality of damages upon these hypersurfaces as proxies for the insurance claims to follow. Finally the model is applied to extensive datasets of hail events, as detected by the NEXRAD network of weather radars.
\end{abstract}


\noindent {\em Journal of Economic Literature Subject Classification:  C46, C12}


\section{Introduction}\label{intro}

\noindent The United States Department of Agriculture (USDA) maintains a large crop insurance program extending to billions of dollars \cite{RMA:A}. Unfortunately, some claims are bound to be fraudulent, and frequently they are related through groups of farmers who act in collusion, extending to conspiring agents, even insurance companies \cite{RLLCS:PC}. Naturally it is desirable to contain this fraud, and there is a need for a good understanding of where actual storm damage has occurred, and to what extent.

To gain a better understanding of hail storm damages, this study investigates the damage to agricultural crops by hail storms, and the pursuant insurance claims. Such claims routinely refer to the distance from the storm center, and are known to respond to countervailing influences. Storm damage occurs with greatest intensity at the center, tapering to insignificance at distance. However, the total of claims filed for damage at the center is small, and increases as more and more claimants reside at greater distances from the center. The total claim value consequently increases from zero as a function of distance to a single mode, and then decreases again to zero.  The research question, therefore, is, ``What model based on fundamentals faithfully replicates this experience?'' The proposed distribution answers this question with parsimony, and is herewith advanced.

This paper is organized as follows. The upcoming section analyzes the log-normal distribution model, which was previously used to describe hail storm damage \cite{L:DMS, LJLRS:CCI}. The following section discusses the 2-dimensional case, under the simplifying assumption that the hail storm does not move over the landscape. The model is that of the independent bivariate normal probability measure of damage intensity. Insofar as damage intensity is independent of direction from the storm center there is only one independent variable --- the radial distance from the center. The resulting marginal distribution on the identity random variable of radius is the $\chi$-distribution. In the next section we extend the model to the traveling form, introducing dependence in the bivariate normal probability measure, and subsequently extend this to the multivariate case. The final phase of the study applies the model to extensive data sets of hail events and their `severe probabilities,' as detected by the NEXRAD network of weather radars.

\section{A log-normal distribution model}\label{sec:lognormal}

\noindent Hail storms can give rise to various forms of damage, including damage to motor vehicles \cite{Schmock99} and to agriculture. In the context of agriculture, the log-normal distribution has been used to describe insurance claim data \cite{L:DMS, LJLRS:CCI}. Although this distribution fits the data reasonably well, we show that there is a theoretical objection to using the log-normal distribution.

The log-normal distribution with parameters $\mu$ and $\sigma$ has density function \cite{AB:LN}
\begin{equation}\label{eq:lognormaldensity}
g_R(r) = g_R(r; \mu, \sigma) = \frac{1}{r \sigma \sqrt{2 \pi}}\, \exp\left(-\frac{(\ln r - \mu)^2}{2\sigma^2} \right),
\qquad r > 0,
\end{equation}
and distribution function
\begin{equation}\label{eq:lognormaldistribution}
G_R(r) = G_R(r; \mu, \sigma) = \NNN\left(\frac{\ln r - \mu}{\sigma}\right),\qquad r > 0,
\end{equation}
with $\NNN$ the standard normal distribution function.

Suppose that $g_R(r)$ is the marginal probability density function in the radial direction of some joint density $g(r,\phi)$ of random variables $R$ and $\Phi$. The damage density at the center can be expressed as the average density over a small disc centered at the center, that is,
\begin{align*}
g(0,\phi) = &\lim_{\varepsilon\to 0} \frac{1}{\pi \varepsilon^2} \int_0^\varepsilon \int_0^{2\pi} g(r,\phi) \dd\phi\dd r \\
 =& \lim_{\varepsilon\to 0} \frac{1}{\pi \varepsilon^2} \int_0^\varepsilon \frac{1}{r \sigma \sqrt{2 \pi}}\, \exp\left(-\frac{(\ln r - \mu)^2}{2\sigma^2} \right) \dd r \\
 =& \lim_{\varepsilon\to 0} \frac{1}{\pi \varepsilon^2} \int_{-\infty}^{(\ln \varepsilon - \mu)/\sigma} \frac{1}{\sqrt{2\pi}} \exp\left( - \frac{1}{2} x^2 \right) \dd x\\
 =& \lim_{\varepsilon\to 0} \frac{1}{\pi \varepsilon^2}
\cdot \NNN\left(\frac{\ln \varepsilon - \mu}{\sigma}\right) \\
 = &\ 0,
\end{align*}
where the last equality follows by applying l'H\^opital's rule. In other words, the log-normal distribution corresponds to a damage pattern with zero damage density in the center, which is unlikely to be the case for a hail storm. This might, however, be desirable for other kinds of storms, like tornados and hurricanes.

\section{A binormal damage pattern and the $\chi$-distribution}\label{sec:binormal}

\noindent If the log-normal distribution is unfit for describing damages, what other distribution is suitable? We make the following desirable assumptions in the damage pattern of a hail storm. The damage function is unimodal at the center, smooth, dependent only on the distance from the center, and scalable to a probability density function.

The simplest distribution with these attributes is the standard bivariate normal, or simply binormal, distribution. We consider the standard probability space $\{\RR^2,\Bmeas,\Pmeas\}$, wherein the first component is the Euclidean plane, the second the Borel sigma algebra, and the third is the binormal independent probability measure. Equip the plane with Cartesian coordinates $(x,y)$ and polar coordinates $(r,\theta)$ and define a random variable $R$ as the identity function on the radial coordinate $r$, independent of $\theta$.  Thus $R$ graphs to an inverted cone with apex at the origin of the plane. One also may define the random variable $\Theta$ as the identity on the angular coordinate $\theta$, independent of $r$. This variable has the uniform distribution.

The usual Euclidian expression of the density of the binormal distribution, founded on the identity random variables $(X,Y)$ on the respective axes with variables $(x,y)$, is
\begin{align*}
f(x,y) &= \frac1{2\pi}\exp\left({-\frac{x^2+y^2}2}\right).
\intertext{The corresponding polar expression is}
g(r) &= \frac1{2\pi}r\exp\left({-\frac{r^2}2}\right).
\end{align*}

Our attention turns to the distribution of the storm damage as distance from the storm center, insofar as the intensity is independent of the direction from the center. The marginal distribution of $R$ in these circumstances is
\[
G(r)
= \Pr\{R\leq r\}
= \frac1{2\pi}\int_0^r s\exp\left({-\frac{s^2}2}\right)\diffns s
= 1-\exp\left({-\frac{r^2}2}\right).
\]
This is the familiar $\chi$-distribution with two degrees of freedom. 

\section{Traveling form of the hail storm damage model}

\noindent Let us assume that at any moment in time, the \emph{damage density} $D_\bfc(\bfx)$ at the location $\bfx \in \RR^2$ of a hail storm is binormally distributed, that is,
\[ D_\bfc(\bfx) = \frac{1}{2\pi} \exp\left(-\frac{1}{2} \|\bfx - \bfc\|^2\right),\]
where $\bfc \in \RR^2$ is the \emph{center} of the hail storm and $\|\cdot\|$ is the Euclidean norm. During the storm, let us assume that the center moves with a constant velocity vector $\bfv\in \RR^2$. Choosing coordinates $\bfx$ such that the center is at the origin $\bf0$ at time $t = 0$, the trajectory of the center is then given by $\bfc = t\bfv$. The \emph{intensity} $I(t)$ of the storm at time $t$ is assumed to be normal,
\[ I(t) = \frac{1}{\sqrt{2\pi}\sigma} \exp\left(-\frac{t^2}{2\sigma^2}\right), \]
with the time coordinate chosen such that the peak intensity happens at time $t = 0$. After scaling the time coordinate by a factor $\sigma$, we can assume that $\sigma = 1$.

Under these assumptions, the \emph{total damage density} $T(\bfx)$ at the point $\bfx$ is given by the marginal density
\[
T(\bfx)
= \int_{-\infty}^\infty I(t) D_{t\bfv}(\bfx) \dd t
= \frac{1}{ (2\pi)^{3/2} } \int_{-\infty}^\infty \exp\left(-\frac{1}{2}\big(t^2 + \|\bfx - t\bfv\|^2\big)\right) \dd t.\\
\]
The integral can be computed by completing the square. Writing
\[ \alpha := \sqrt{1 + \|\bfv\|^2},\qquad s := \alpha t - \frac{\langle\bfv, \bfx\rangle}{\alpha}, \]
with $\langle\cdot,\cdot\rangle$ the standard inner product, one finds that the total damage
\begin{align*}
T(\bfx) & = \frac{1}{2\pi\alpha} \exp\left(-\frac12 \left(\|\bfx\|^2 - \frac{\langle\bfv,\bfx\rangle^2}{\alpha^2} \right)\right)\cdot \frac{1}{\sqrt{2\pi}} \int_{-\infty}^\infty \exp\left(-\frac12 s^2\right)\dd s \\
        & = \frac{\exp\left(-\frac12 \left(\|\bfx\|^2 - \frac{\langle\bfv,\bfx\rangle^2}{1 + \|\bfv\|^2} \right)\right)}{2\pi\sqrt{1 + \|\bfv\|^2}}
\end{align*}
is also binormally distributed, but now with a correlation in its random vector. To bring this density in standard form, write $\bfv = (v_1, v_2)$ and introduce the parameters
\[
\sigma_1 := \sqrt{1 + v_1^2},\qquad
\sigma_2 := \sqrt{1 + v_2^2},\qquad
\rho := \frac{v_1 v_2}{\sqrt{(1 + v_1^2)(1 + v_2^2)}}.
\]
Then
\[ T(\bfx) = \frac{1}{2\pi\sqrt{\det(\bfSigma)}} \exp\left(-\frac{1}{2}\bfx^\TT \bfSigma^{-1}\bfx \right), \]
which is the standard form of the bivariate normal distribution with zero mean and covariance matrix
\[
\bfSigma = 
\begin{bmatrix}
            \sigma_1^2 & \rho \sigma_1 \sigma_2  \\
\rho\sigma_1\sigma_2   &               \sigma_2^2
\end{bmatrix}.
\]

\section{The marginal distribution at a distance}

\noindent As for the standard binormal distribution in Section \ref{sec:binormal}, we wish to reduce the dimension of the damage density $T(\bfx)$, so as to cancel out fluctuations in the data and arrive at a more regular and simpler distribution. The damage density, with general covariance matrix $\bfSigma$, is no longer constant in the angular direction. It is therefore natural to change to a coordinate system in which the damage density becomes standard binormal, before reducing the dimension. We choose to perform this process for a general multinormal distribution, as this is not much harder than the bivariate case and might be used for other modeling purposes.

Suppose a random vector $\bfX$ on $\RR^n$ is multivariate normally distributed with density function
\begin{equation}\label{eq:multinormaldistribution}
T_\bfX(\bfx) = \frac{1}{(2\pi)^{n/2}\sqrt{\det(\bfSigma)}}\exp\left(-\frac{1}{2}(\bfx-\bfmu)^\TT\bfSigma^{-1} (\bfx-\bfmu) \right),
\end{equation}
with mean $\bfmu$ and covariance matrix $\bfSigma$ of full rank $n$. If we are given data $(P_i, \bfx_i)$, where $P_i$ is the frequency of the event $\bfx_i$, one has the maximum likelihood estimators
\begin{equation}\label{eq:EstimateMean}
\widehat{\bfmu}
:= \frac{1}{\sum_i P_i} \sum_i P_i \bfx_i,
\end{equation}
\begin{equation}\label{eq:EstimateCovarianceMatrix}
\widehat{\bfSigma} =\begin{bmatrix} \widehat{\sigma}_x^2 & \widehat{\sigma}_{xy}\\ \widehat{\sigma}_{xy}& \widehat{\sigma}_y^2\end{bmatrix}
:= \frac{1}{\sum_i P_i}
\sum_i P_i
(\bfx_i - \widehat{\bfmu})(\bfx_i - \widehat{\bfmu})^\TT
\end{equation}
for the parameters $\bfmu$, $\bfSigma$.

While the independent binormal distribution has angular symmetry with respect to its center, this is no longer the case for the general multinormal distribution. To reduce to a univariate distribution, we introduce a distance function that is zero at the center and constant along the level curves of $T_\bfX(\bfx)$. It is easily checked that the function $d:\RR^2 \longrightarrow [0,\infty)$ defined by
\begin{equation}\label{eq:distance}
d(\bfx) = \sqrt{ (\bfx - \bfmu)^\TT \bfSigma^{-1} (\bfx - \bfmu) }
\end{equation}
has these properties. In the case of the standard binormal distribution with $\bfmu = \bf0$ and $\bfSigma$ the identity matrix, this is the ordinary Euclidean distance to the origin. We can then consider a marginal distribution in the ``radial distance direction'' $d$.

To derive the density $T_R(r)$ of this marginal distribution, we can change coordinates $\bfx$ such that $T_\bfX$ becomes independent multinormal. The density \eqref{eq:multinormaldistribution} partitions Euclidean $n$-space into level hypersurfaces with constant probability (i.e., with constant distance $d$ from the center). Since the covariance matrix $\bfSigma$ is symmetric positive definite, it admits an orthogonal diagonalization
\[ \bfSigma = \bfQ^\TT \bfD \bfQ, \qquad \bfQ := [\bfv_1, \ldots, \bfv_n], \qquad \bfD := \diag\{a_1^2,\ldots, a_n^2\}. \]
The level hypersurfaces form a family of hyperellipsoids with center $\bfmu$, semi-axis lengths $a_1, \ldots, a_n$ in constant proportion $[a_1:\cdots:a_n]$, and directions of the principal axes given as corresponding eigenvectors $\bfv_1, \ldots, \bfv_n$ of the covariance matrix $\bfSigma$. Transforming to a random vector $\bfY$ by the change of coordinates $\bfy := \bfQ(\bfx-\bfmu)$ yields the probability density function
\begin{eqnarray*}
T_\bfY(\bfy) & = & T_\bfX(\bfQ^\TT \bfy) \cdot \left| \det \bfQ^\TT \right| \\
       & = & \frac{1}{(2\pi)^{n/2} a_1 \cdots a_n} \exp \left(-\frac{1}{2}\left(\frac{y_1^2}{a_1^2} + \cdots + \frac{y_n^2}{a_n^2} \right) \right)
\end{eqnarray*}
of $\bfY$. Changing to hyperspherical coordinates by the map
\[ (0,\infty)\times [0,\pi]^{n-2} \times [0,2\pi) \longrightarrow \RR^n\]
defined by
\[
\bfz = \begin{bmatrix}r\\ \phi_1\\ \vdots\\ \phi_{n-2}\\ \phi_{n-1}\end{bmatrix}\longmapsto 
\bfy = 
\begin{bmatrix} y_1\\ y_2\\ \vdots\\ y_{n-1}\\ y_n \end{bmatrix} =
r
\begin{bmatrix}
a_1 \cos(\phi_1)\\
a_2 \sin(\phi_1)\cos(\phi_2)\\ 
\vdots\\
a_{n-1} \sin(\phi_1)\cdots \sin(\phi_{n-2}) \cos(\phi_{n-1})\\
\ \ \ a_n \sin(\phi_1)\cdots \sin(\phi_{n-2}) \sin(\phi_{n-1})
\end{bmatrix}
\]
yields a random vector $\bfZ$ with probability density function
\begin{eqnarray*}
T_\bfZ(\bfz) & = & T_\bfY(\bfy) \cdot \left|\det\left\{\frac{\partial(y_1,y_2,\ldots,y_n)}{\partial(r,\phi_1,\ldots,\phi_{n-1})}\right\} \right|\\
       & = & \frac{1}{(2\pi)^{n/2}} r^{n-1}\sin^{n-2}(\phi_1) \sin^{n-3}(\phi_2) \cdots \sin(\phi_{n-2})\exp \left(-\frac{r^2}{2}\right)
\end{eqnarray*}
that respects the foliation by hyperellipsoids. Note that
\[ r^2 = \frac{y_1^2}{a_1^2} + \cdots \frac{y_n^2}{a_n^2} = \bfy^\TT \bfD^{-1} \bfy = (\bfx - \bfmu)^\TT \bfQ^\TT \bfD^{-1} \bfQ(\bfx - \bfmu) = d(\bfx)^2, \]
meaning that the distance $d$ in $\bfx$-space is the radius $r$ in $\bfz$-space.

Integrating over the angular random variables, one is left with the marginal radial random variable $R$ with marginal probability density function 
\begin{eqnarray*}
T_R(r) & = & \int_0^\pi \cdots \int_0^\pi \int_0^{2\pi} T_\bfZ(r,\phi_1,\ldots,\phi_{n-2},\phi_{n-1})\dd \phi_{n-1} \dd \phi_{n-2} \cdots \dd \phi_1\\
       & = & \frac{2\pi^{n/2}}{\Gamma\big(\frac{n}{2}\big)}\cdot \frac{1}{(2\pi)^{n/2}} r^{n-1}\exp \left(-\frac{r^2}{2}\right) =
             \frac{2^{1-n/2}}{\Gamma\big(\frac{n}{2}\big)} r^{n-1}\exp \left(-\frac{r^2}{2}\right),
\end{eqnarray*} 
where we used that the surface area of the unit sphere of dimension $n~-~1$ is $2\pi^{n/2}/\Gamma(n/2)$, with $\Gamma$ the gamma function. One recognizes $T_R(r)$ as the density function of the $\chi$-distribution with $n$ degrees of freedom. One hits the interior of the hyperellipsoid defined by $\frac{y_1^2}{a_1^2} + \cdots + \frac{y_n^2}{a_n^2} = R^2$ with probability
\[ \Pr(0\leq R\leq r) = \int_0^r T_R(s)\dd s = P(n/2, r^2/2), \]
where $P$ is the regularized Gamma function \cite[\S 6.5.1]{AS:HB}.

For $n=2$ we recover the hail storm setting. To evaluate insurance claims it is helpful to compare, at the point $\bfx$, the reported total damage to the expected total damage. Since the latter quality is, by definition, constant along the level curve through $\bfx$, it is tempting to reduce the dimension of the problem by considering the marginal distribution in the radial direction, which has density function
\[ T_R(r) = r\exp \left(-\frac{r^2}{2}\right) \] 
corresponding to the $\chi$-distribution with two degrees of freedom, also known as the Rayleigh distribution. The total damage within the ellipse defined by $\frac{y_1^2}{a_1^2} + \frac{y_2^2}{a_2^2} = R^2$ takes on the particularly simple form
\[ \Pr(0\leq R\leq r) = \int_0^r T_R(s)\dd s = 1 - \exp\left(-\frac{r^2}{2}\right). \]

Since the $\chi$-distribution is the marginal distribution in the radial direction of a bivariate distribution with mode at its center, it does not suffer from the theoretical objection to the log-normal distribution raised in Section \ref{sec:lognormal}. As we shall see in the example data in the next section, even though the log-normal distribution depends on an additional parameter, it only has a somewhat better overall fit than the $\chi$-distribution. In addition its density is too low near the center and its tail is too fat.

\begin{table}[h!]
{\small
\begin{tabular}{ccccccccc}
\toprule
$t_i$ & $\bfx_i^\TT$ & $P_i$ & & $t_i$ & $\bfx_i^\TT$ & $P_i$\\
\midrule
16:56:04 & [-89.72179, 31.44091] & 0.6 & & 17:46:33 & [-89.38189, 31.64263] & 0.8\\
17:00:05 & [-89.68649, 31.46652] & 0.4 & & 17:47:19 & [-89.39462, 31.62080] & 0.2\\
17:00:16 & [-89.69341, 31.46517] & 0.8 & & 17:50:45 & [-89.35648, 31.63880] & 0.5\\
17:00:40 & [-89.67842, 31.47161] & 0.2 & & 18:07:34 & [-89.22546, 31.65128] & 0.8\\
17:04:29 & [-89.62775, 31.48436] & 0.7 & & 18:11:46 & [-89.18574, 31.66288] & 0.5\\
17:05:04 & [-89.62963, 31.47485] & 0.1 & & 18:24:23 & [-89.03573, 31.68747] & 0.5\\
17:05:52 & [-89.63511, 31.49470] & 0.3 & & 18:27:45 & [-89.04506, 31.64438] & 0.1\\
17:08:41 & [-89.57959, 31.50128] & 0.7 & & 18:28:35 & [-88.98864, 31.65766] & 0.5\\
17:12:38 & [-89.55715, 31.51726] & 0.3 & & 18:28:35 & [-88.94628, 31.84381] & 0.1\\
17:12:54 & [-89.61882, 31.50527] & 0.7 & & 18:32:48 & [-88.91338, 31.78937] & 0.5\\
17:17:06 & [-89.57194, 31.52256] & 0.9 & & 18:32:48 & [-88.97600, 31.67251] & 0.3\\
17:18:25 & [-89.55338, 31.53365] & 0.4 & & 18:37:00 & [-88.86683, 31.83042] & 0.4\\
17:21:18 & [-89.49632, 31.55593] & 0.8 & & 18:45:32 & [-88.74685, 31.87335] & 0.6\\
17:22:42 & [-89.51342, 31.55574] & 0.2 & & 18:58:05 & [-88.75892, 31.71739] & 0.5\\
17:24:12 & [-89.49149, 31.57349] & 0.3 & & 19:02:18 & [-88.71288, 31.72012] & 0.4\\
17:25:31 & [-89.46705, 31.57082] & 0.7 & & 19:10:01 & [-88.48747, 31.96443] & 0.4\\
17:27:07 & [-89.46227, 31.58279] & 0.1 & & 19:10:44 & [-88.49454, 31.97819] & 0.3\\
17:29:44 & [-89.43840, 31.58658] & 0.6 & & 19:12:36 & [-88.47559, 31.98316] & 0.5\\
17:33:57 & [-89.42432, 31.59477] & 0.6 & & 19:14:54 & [-88.44588, 31.99579] & 0.5\\
17:35:45 & [-89.46156, 31.58444] & 0.1 & & 19:14:57 & [-88.45017, 31.99291] & 0.3\\
17:35:56 & [-89.46227, 31.58279] & 0.1 & & 19:17:11 & [-88.40633, 31.99112] & 0.6\\
17:38:09 & [-89.42432, 31.59477] & 0.6 & & 19:19:09 & [-88.38118, 32.02913] & 0.3\\
17:42:21 & [-89.42432, 31.59477] & 0.9 & & 19:19:47 & [-88.40560, 32.02777] & 0.4\\
\bottomrule
\end{tabular}
}
\caption{Hail events belonging to a hail storm on January 20, 2010, in the vicinity of Laurel, Mississippi. Each of the 46 hail events lists a time $t_i$, a location $\bfx_i$ as a column vector $[$longitude, latitude$]^\TT$, and a severe probability~$P_i$.}\label{tab:haildata}
\end{table}

\section{Fitting the model to data}\label{sec:ModelFit}

\begin{figure}
\begin{center}
\begin{overpic}[grid=False, width=\columnwidth, height=0.5\columnwidth]{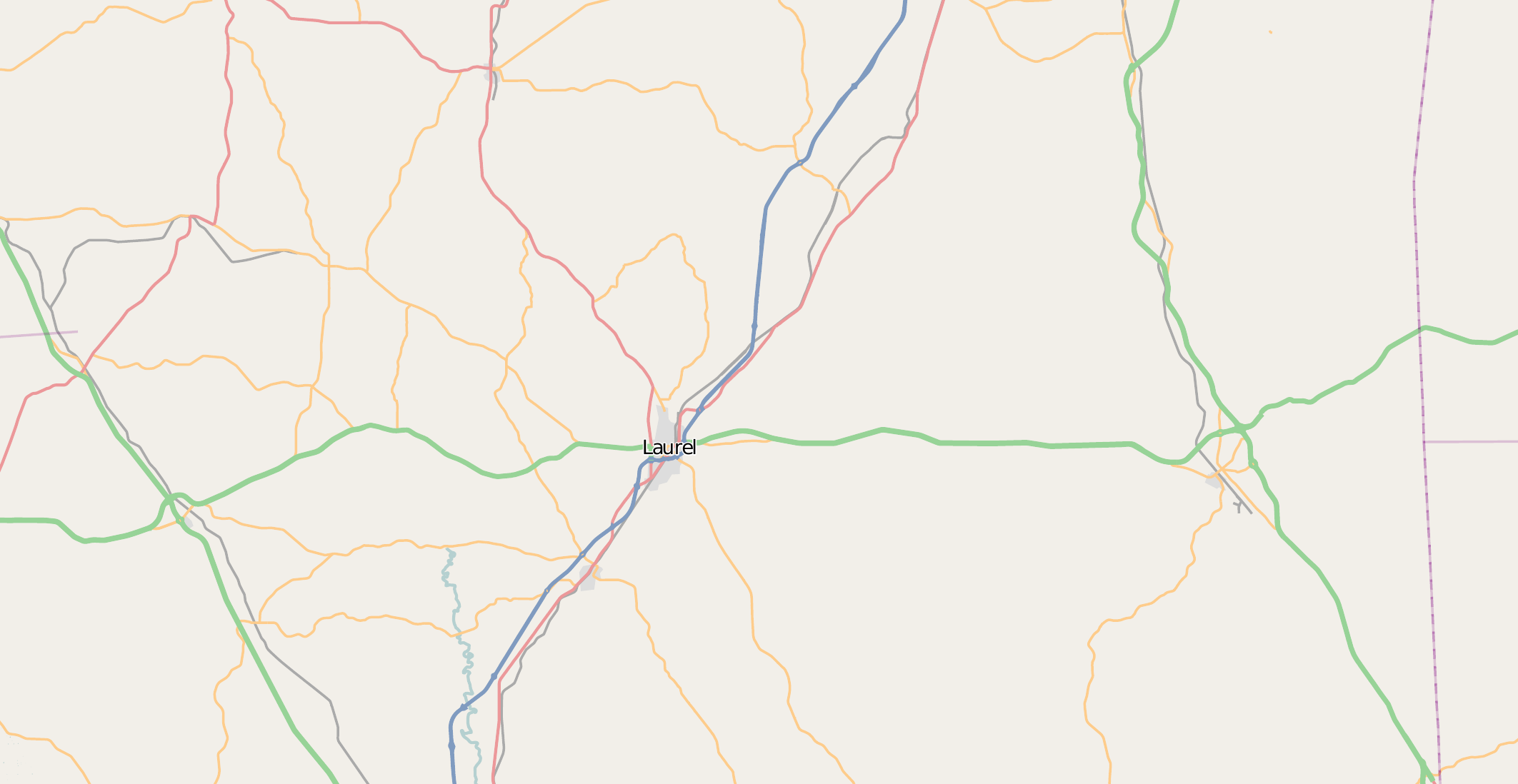}
\includegraphics[bb=15 10 405 182, width=\columnwidth, height=0.5\columnwidth]{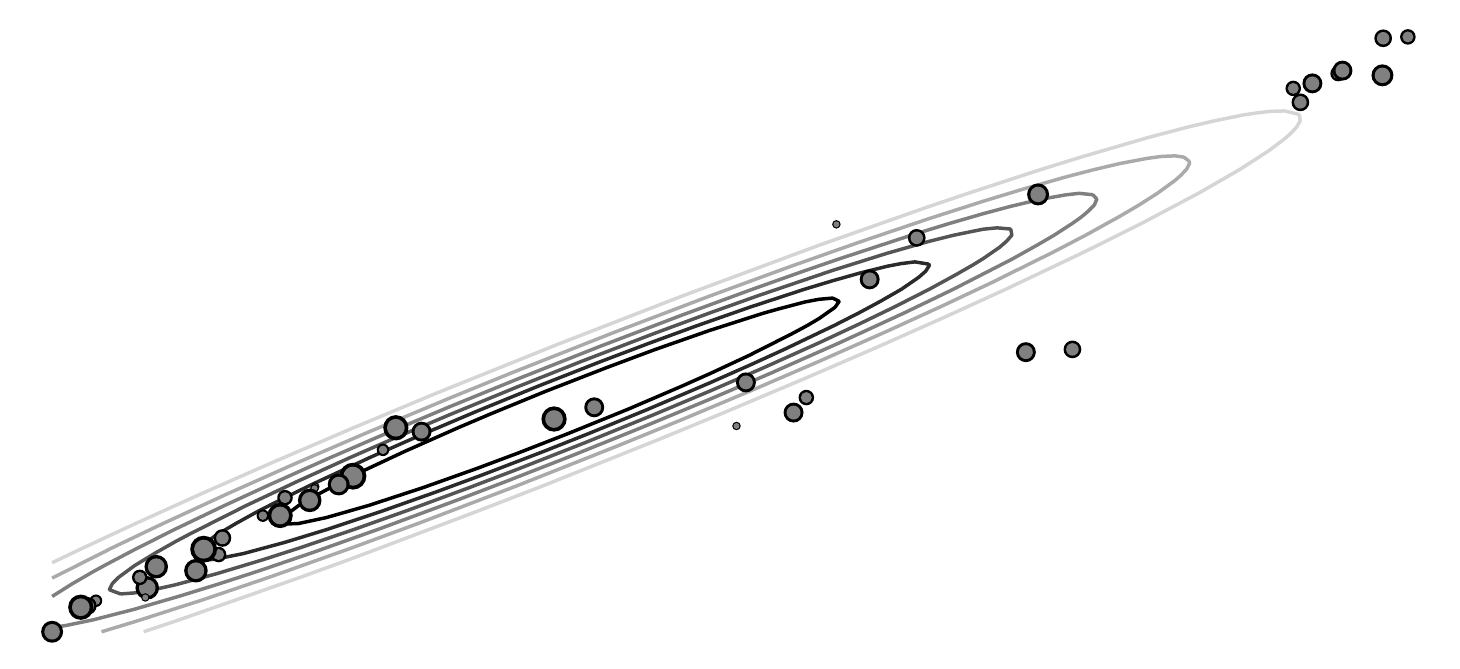}
\end{overpic}
\end{center}
\caption{Drawn as points on top of a map in the vicinity of Laurel, Mississippi, are the hail events from Table \ref{tab:haildata}, with sizes proportional to their severe probabilities. In addition, contour lines of a fitted binormal distribution are drawn.}\label{fig:OpenStreetMap}
\end{figure}

\noindent In this section we fit our model to a data set of hail events, as estimated by the Next Generation Weather Radar system (NEXRAD) network \cite{NEXRAD}. Distributed throughout the United States and selected overseas locations, over a hundred weather radars measure the reflectivity, mean radial velocity, and spectrum width. These meteorological base data quantities are used to search for patterns that estimate the presence, and likelihood, of various kinds of severe weather events. One of the data sets derived from this processing is the Hail Index Overlay, which is designed to locate storms with the potential to produce hail. This data set is organized as a collection of hail events and the probability that the event is severe, which can be thought of as a potential intensity of the hail event. The National Climatic Data Center makes these hail events publicly available through the Severe Weather Data Inventory \cite{NCDC:SWDI}.

We are, however, not interested in single hail events, but in hail storms. We have chosen the single-linkage distance hierarchical agglomerative clustering method to define what is meant by a single hail storm, as we found this to be in line with our own intuitive notion of a storm. Using \texttt{R} \cite{R} and in particular the package \texttt{flashClust} \cite{Langfelder.Horvat12}, we compute the hierarchical clustering tree from a large collection of hail events in January, 2010. See Murtagh \cite{Murtagh83} for the details of the underlying algorithm. \emph{A priori} we do not know how many storms to expect. Following a rule of thumb, we cut the dendrogram when the next merging gives rise to a disproportionate jump in the clustering criterion. In this manner, we clustered the hail events in the month January in several storms. We chose one representative storm that was not too large, from January 20, 2010, which is listed in Table \ref{tab:haildata} and shown in Figure \ref{fig:OpenStreetMap} on top of a map of the vicinity of Laurel, Mississippi \cite{OSM:OSM}.

The events appear relatively near each other and far from either pole, implying that we can approximately treat the longitude and latitude as Cartesian coordinates. Let us assume that the locations $\bfx_i$ of the hail events in Table \ref{tab:haildata} are sampled from a binormal distribution with density as in Equation \eqref{eq:multinormaldistribution}. Each $\bfx_i$ comes with a ``severe probability'' $P_i$ that is interpreted as a weight of the event. The center $\bfmu$ and covariance matrix $\bfSigma$ of the storm can be estimated by maximum likelihood, as in \eqref{eq:EstimateMean} and \eqref{eq:EstimateCovarianceMatrix}. The resulting fitted binormal distribution is depicted in Figure \ref{fig:OpenStreetMap} by some of its contour lines.

By the discussion of the previous section, the marginal distribution in the radial direction is the $\chi$-distribution with two degrees of freedom. Because the pair $(a_1, a_2)$ of semi-axes is only defined up to multiplication by a constant, we can consider a family of $\chi$-like distributions
\[ F(r; \lambda) = 1 - \exp\left( -\frac{1}{2} \lambda^2 r^2 \right),\qquad r > 0, \]
parametrized by $\lambda > 0$.

To find the estimator $\widehat{\lambda}$ of the parameter $\lambda$ that fits our data best, we reorder the data by distance $d$ in \eqref{eq:distance} from the center. Let $\pi$ be a permutation of the indices of the hail events for which $\big(d(\bfx_{\pi_i}) \big)_i$ becomes a nondecreasing sequence of distances. Estimating the parameter $\lambda = \widehat{\lambda}$ for which $F(r; \lambda)$ is the best fit of our data can be done by solving the nonlinear least square problem
\begin{equation}\label{eq:LeastSquaresChi}
\widehat{\lambda} := \underset{\lambda > 0}{\operatorname{argmin}}\ \sum_i \Big[F\big(d(\bfx_{\pi_i}); \lambda \big) - \sum_{\pi_j \leq \pi_i} P_{\pi_j} \Big]^2.
\end{equation}
Solving this problem numerically using \texttt{Sage} \cite{Sage}, we find that the sum of squares reaches its minimum of $S_F = 0.067$ at $\widehat{\lambda} \approx 7.308$.

To compare the $\chi$-distribution to the log-normal distribution, we can either fit the log-normal distribution using the Euclidean distance or using our distance function $d$ estimated in \eqref{eq:distance}. For instance, in the latter case a best fitting log-normal distribution can be found by numerically solving the nonlinear least square problem
\begin{equation}\label{eq:LeastSquaresLogNormal}
(\widehat{\mu}, \widehat{\sigma}) := \underset{(\mu,\sigma) \in \RR\times (0,\infty)}{\operatorname{argmin}}\ \sum_i \Big[G\big(d(\bfx_{\pi_i}); \mu, \sigma\big) - \sum_{\pi_j \leq \pi_i} P_{\pi_j} \Big]^2.
\end{equation}
One finds that the sum of squares reaches its minimum of $S^d_G = 0.0483$ at $\widehat{\mu} \approx -1.862$ and $\widehat{\sigma} \approx 0.6227$. For the hail storms recorded in January 2010, penalties of a fitted $\chi$- and log-normal distribution can be found in Table~\ref{tab:penalties}. In this table, storm 4 refers to the storm analyzed in detail in this paper.

\begin{table}
\begin{center}
\begin{tabular*}{\columnwidth}{c@{\extracolsep{\stretch{1}}}*{10}{r}@{}}
\toprule
Storm     &    1 &    2 & 3 & 4 & 5 & 6 & 7 & 8 & 9 & 10\\
\midrule
Events    &   42 &   54 &   45 &    46 & 3426 & 1022 &   36 &  192 &   58 &   76\\
$S_F$     & 0.25 & 0.43 & 0.21 & 0.067 & 23.7 &  3.2 & 0.16 & 0.23 & 0.41 & 0.15\\
$S_G$     & 0.17 & 0.23 & 0.17 & 0.045 &  2.5 &  1.9 & 0.08 & 0.17 & 0.05 & 0.24\\
$S^d_G$   & 0.09 & 0.29 & 0.20 & 0.048 &  4.5 &  2.2 & 0.07 & 0.08 & 0.05 & 0.33\\
\bottomrule
\end{tabular*}
\end{center}
\caption{For the ten hail storms recorded in January 2010, number of events and penalties of fitted $\chi$-distributions ($S_F$) and log-normal distributions using the Euclidean distance function ($S_G$) and the distance function $d$ ($S^d_G$).}\label{tab:penalties}
\end{table}

Comparing the sums of squares $S_G$ and $S_F$, the log-normal distribution has a somewhat better overall fit than the $\chi$-distribution, which is to be expected because of its additional parameter. Plotting the residuals of the fitted $\chi$-distribution and log-normal distribution shows that they are approximately normally distributed. The F-test of the equality of two variances yields an F-statistic of approximately $0.067/0.0483$ with corresponding P-value $0.142$, taking into account the additional parameter of the log-normal distribution. The null-hypothesis of equality of variance can therefore not be rejected at the 10\% significance level. 

\begin{figure}
\begin{center}
\includegraphics[scale=0.7]{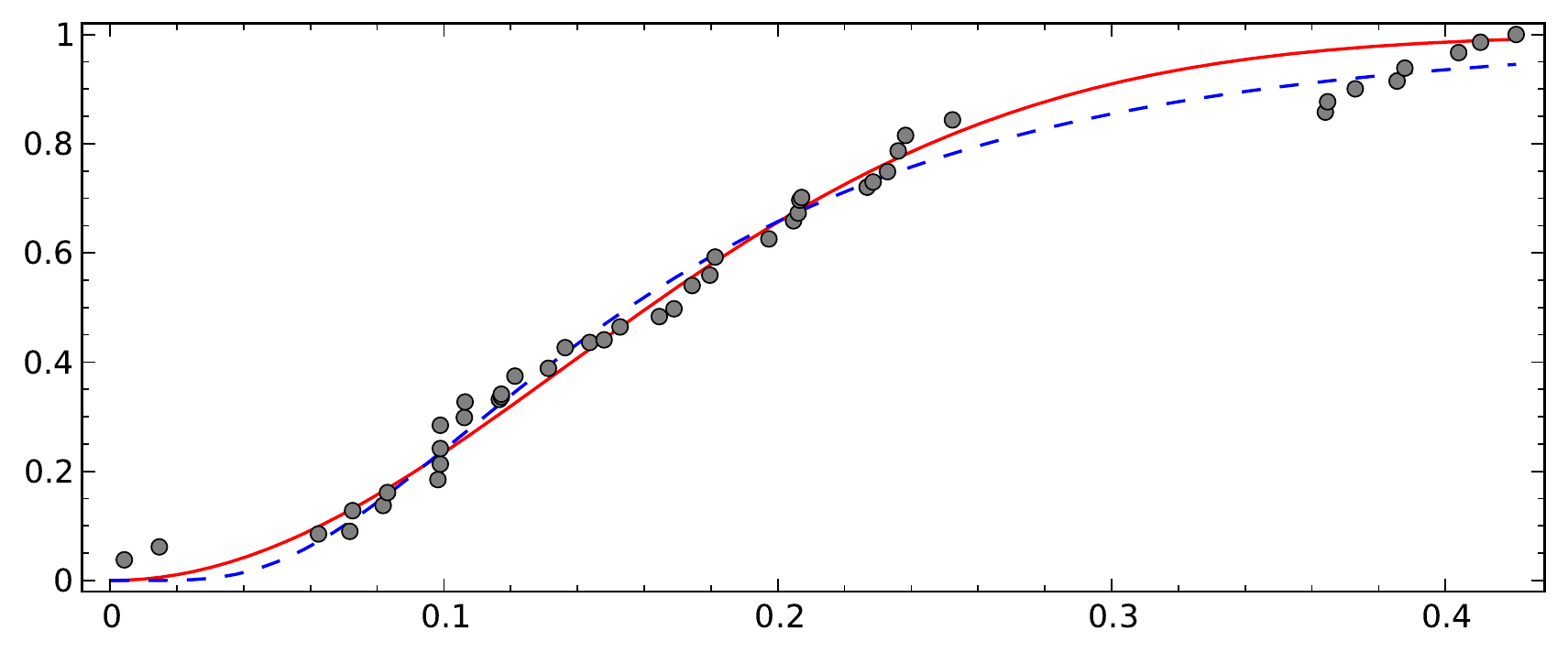}
\end{center}
\caption{The empirical distribution for the distance function \eqref{eq:distance}, together with a fitted $\chi$-distribution (drawn solid) and a fitted log-normal distribution (drawn dashed), found by solving the nonlinear least squares problems \eqref{eq:LeastSquaresChi} and \eqref{eq:LeastSquaresLogNormal}.} \label{fig:FittedChiDistributionLogNormal}
\end{figure}

\begin{figure}
\begin{center}
\includegraphics[scale=0.7]{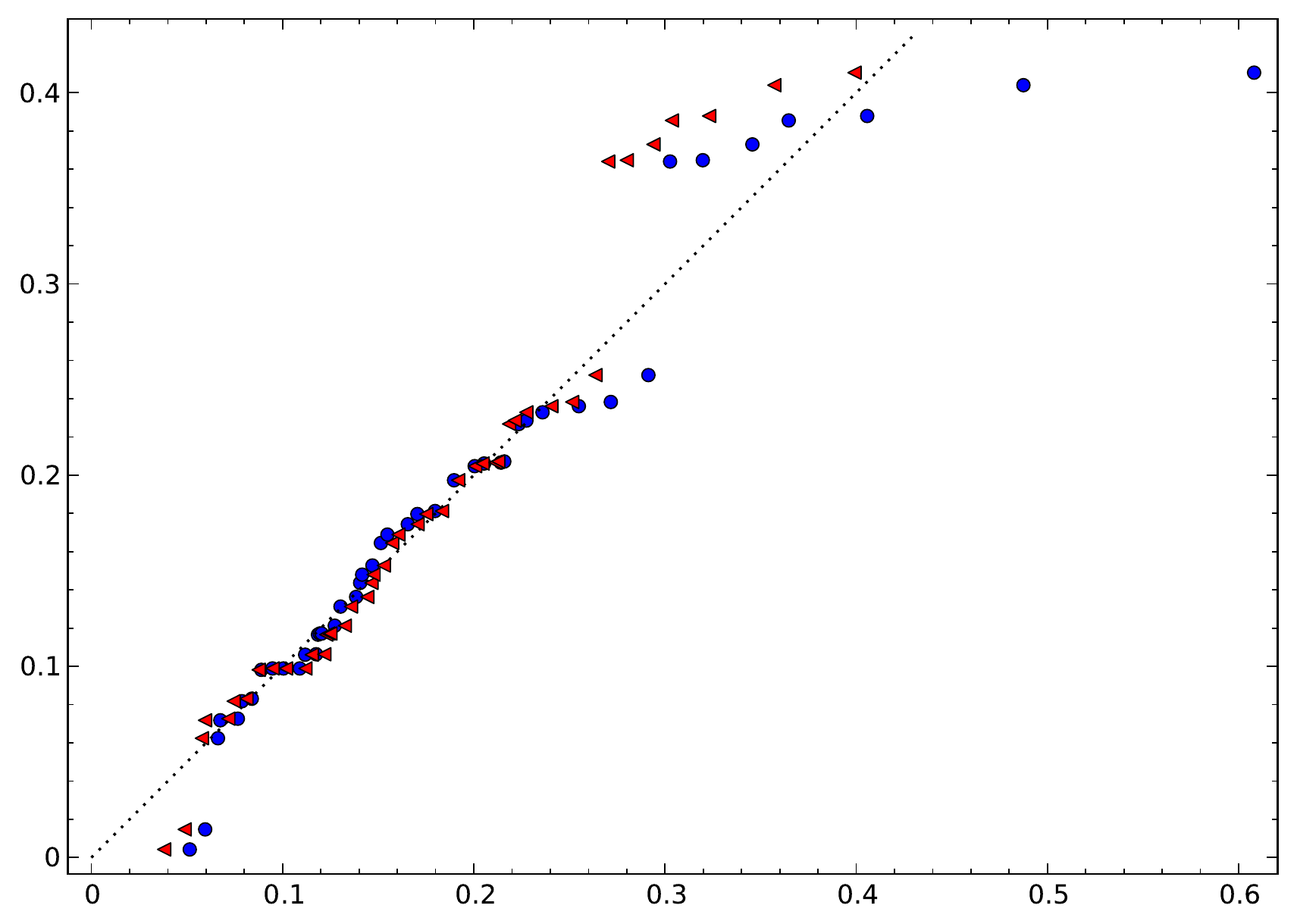}
\end{center}
\caption{A Q--Q plot comparing the empirical distribution on the vertical axis to the fitted $\chi$-distribution (${\color{red} \blacktriangleleft}$) and the fitted log-normal distribution (${\color{blue} \bullet}$) on the horizontal axis.}\label{fig:QQplot}
\end{figure}

Figure \ref{fig:FittedChiDistributionLogNormal} simultaneously shows the empirical distribution for the distance function \eqref{eq:distance}, the best-fitted $\chi$-distribution and best-fitted log-normal distribution. Qualitatively, the fitted log-normal distribution is too low near the origin, confirming the discussion in Section \ref{sec:lognormal}, and its tail seems to be too fat for the data. This can be seen more clearly from the Q--Q plot in Figure~\ref{fig:QQplot}. Note that the $\chi$-distribution is also too low near the origin, but somewhat better than the log-normal distribution. The log-normal distribution tends to provide a better fit in the middle.

Finally, let us note some limitations of the model. In order to approximate longitude and latitude by Cartesian coordinates, the storm cannot be too large. In addition, for the Coriolis effect to be negligible, the storm cannot last too long. When using these hail intensities as proxies for damage claims, the underlying topography should be homogeneous. This is for instance the case with large-scale corn field agriculture. Moreover, our model does not reflect that different types of hail storms can cause different types of damage \cite{Schiesser97}. For instance, larger hail stones are more likely to damage motor vehicles, while hail storms with small but numerous hail stones have a greater damaging effect on crops.

\section*{Acknowledgments}

\noindent The authors wish to thank Bertis B. Little of Texas A\&{}M University, Tarleton State University Campus and his colleagues, for inspiration on this topic \cite{L:DMS} and fundamental motivating works \cite{LJLRS:CCI,RLLCS:PC}. In addition we express our gratitude to Juan Gerardo Alc\'azar for his help with clustering methods, and Bj\o rn Sundt for his detailed feedback on an earlier version of the paper.

\section*{References}

\bibliographystyle{amsxport}

\end{document}